# ASYMPTOTIC BEHAVIOR OF A METAPOPULATION MODEL

By A. D. Barbour[1] and A. Pugliese[2]

*Universität Zürich and Universitá di Trento*

We study the behavior of an infinite system of ordinary differential equations modeling the dynamics of a metapopulation, a set of (discrete) populations subject to local catastrophes and connected via migration under a mean field rule; the local population dynamics follow a generalized logistic law. We find a threshold below which all the solutions tend to total extinction of the metapopulation, which is then the only equilibrium; above the threshold, there exists a unique equilibrium with positive population, which, under an additional assumption, is globally attractive. The proofs employ tools from the theories of Markov processes and of dynamical systems.

**1. Introduction.** The simplest models of population growth and regulation are formulated in terms of a more or less isolated population in a single habitat. However, the importance of the spatial dimension has been recognized in a number of ecological processes, resulting in one of the most active topics in theoretical ecology: see, for instance, the two recent collections [26] and [8] and the review article by Neuhauser [22]. These ideas have stimulated the development of spatially structured stochastic populations models, as in [23] and [10], whose mathematical analysis is generally very hard.

A very simple model recognizing the spatial dimension of ecological processes was introduced by Levins [19] in 1969. He envisaged a metapopulation consisting of many distinct habitat patches, within each of which the population behaves much as in the single population models, but which are linked to one another by migration. In his highly simplified model, patches are designated as occupied or not, and all occupied patches are taken to

Received December 2003; revised July 2004.
[1]Supported in part by Schweizer Nationalfonds Projekt 20-61753.00.
[2]Supported in part by CNR Grant 00.0142.ST74 (Italy).

be equivalent, irrespective of the number of individuals present. With these simplifications, he obtained a single differential equation,

$$\frac{dp}{dt} = c_L p(1-p) - \nu_L p, \tag{1.1}$$

describing the behavior of the system: here, $p = p(t)$ represents the proportion of occupied patches, $\nu_L$ is the extinction rate and $c_L$ is the colonization rate per occupied patch. Hence, an equilibrium exists only if $c_L > \nu_L$, and, in that case, the proportion of empty patches at equilibrium is $\nu_L/c_L$. His ideas have been widely used, both in theoretical papers and in wildlife management problems (see, e.g., [15]).

Levins' metapopulation model has two major weaknesses: on the one hand, it is based on a mean field assumption (the colonization rate in a patch depends only on the overall proportion of patches occupied); on the other hand, all patches are assumed to be equal and described simply as empty or occupied, disregarding local population dynamics. Addressing the first issue requires the consideration of spatial stochastic processes as mentioned above. For the second, some authors have generalized Levins' model by taking into account the numbers of individuals in the occupied patches, giving rise to the so-called structured metapopulation models [12]: they consist either of a finite [21] or infinite number of ordinary differential equations [5], or of a partial differential equation [12, 13], where the structuring variable $x$ represents the number of individuals per patch. However, very few analytical results are available for models of complexity comparable to ours, and the behavior of these models has mainly been explored through simulation.

In this paper we investigate the deterministic approximation to the metapopulation model discussed in [1]. This is a stochastic mean field metapopulation model, in which the number of individuals in a patch is governed by a birth, death and catastrophe process, with the same transition rates in each patch, together with migration between the patches with a uniform transition rate $\gamma$ per individual, destinations being chosen uniformly at random among all patches. This last, mean field assumption is probably the least biologically realistic, but has been used in several papers [15], and may make very good sense for metapopulations of parasites in which the patches represent host animals. At all events, it makes the mathematical treatment substantially simpler.

As is shown in [1], when the number of patches becomes very large, one can approximate the stochastic model with the following infinite system of differential equations:

$$p_i'(t) = -\left[(b_i + d_i + \gamma)i + \nu + \rho\gamma \sum_{j=0}^{\infty} j p_j(t)\right] p_i(t)$$

$$
\begin{aligned}
&+ \left[ b_{i-1}(i-1) + \rho\gamma \sum_{j=0}^{\infty} j p_j(t) \right] p_{i-1}(t) \\
&+ [d_{i+1} + \gamma](i+1)p_{i+1}(t), \qquad i \geq 1,
\end{aligned}
$$

(1.2)

$$p_0'(t) = \nu\left( \sum_{j=0}^{\infty} p_j(t) - p_0(t) \right) + (d_1 + \gamma)p_1(t) - \rho\gamma\left( \sum_{j=0}^{\infty} j p_j(t) \right) p_0(t),$$

$$\underline{p}(0) = \underline{p}^0,$$

in which $p_i(t)$ denotes the proportion of patches that are occupied by $i$ individuals, $i \geq 0$. The parameters $b_i$ and $d_i$ represent the per capita birth and death rates in a patch occupied by $i$ individuals, the catastrophe rate is $\nu$ in each patch, the migration rate is $\gamma$ per individual, and $\rho$ is the probability of a migrant successfully reaching another patch. Note that this model is very similar to those studied by Metz and Gyllenberg [21] as structured metapopulation models with finite patch size, and by Casagrandi and Gatto [5].

We also assume the following:

(H1) $i b_i$ is concave and nondecreasing; $i d_i$ is convex and nondecreasing.

It can easily be seen that (H1) implies that $b_i$ is nonincreasing and $d_i$ nondecreasing. Hence, there exist $b_\infty = \lim_{i \to \infty} b_i$ and $d_\infty = \lim_{i \to \infty} d_i$, for which we further assume that

(H2) $b_\infty < d_\infty + \gamma(1-\rho) + \nu$.

Generally, in logistic demography, the existence of a carrying capacity is assumed: that is, there is a value $K$ such that $b_K = d_K$, which automatically implies that $b_\infty < d_\infty$. (H2) is weaker than that, and is, in fact, the natural condition: if $b_\infty \geq d_\infty + \gamma(1-\rho) + \nu$, there can be no nontrivial equilibrium, as is proved in Proposition 3.4.

The assumptions of concavity of $i b_i$ and convexity of $i d_i$ are satisfied in many examples, but not in all; for instance, a Ricker-type birth function $b_i = b_0 \exp\{-\beta i\}$ is not allowed. However, they are mathematically convenient assumptions, if the uniqueness of any nontrivial equilibrium solution to equations (1.2) is to be guaranteed, and we make use of them in several steps of our proofs; they could certainly be relaxed, but it is not easy to see what general conditions would better replace them.

The existence and uniqueness of the solutions to (1.2) are established in [1], and a summary of those of her results relevant to this paper is given at the end of the section. In this paper we consider the possible equilibria $\pi$ of (1.2), using stochastic coupling arguments that are developed in Section 2. There is always the "extinction" equilibrium, with $\pi(0) = 1$ and $\pi(i) = 0$, $i \geq 1$; this is also the eventual limit of all finite patch stochastic systems, and makes the theory of quasi-equilibria of essential importance for such models. In addition, if a threshold condition is satisfied, we show in Section 3 that there is a unique nonnegative equilibrium having $\pi(0) < 1$ (Theorem 3.1). This distribution is shown to be the equilibrium distribution for the single patch dynamics, in which immigration from outside is fixed at a constant "effective" rate, determined by the nonzero solution of a fixed point equation (3.3). In Theorem 4.5 of Section 4, we prove *global* convergence to this equilibrium when the threshold condition is satisfied, under the additional assumption that $d_\infty < +\infty$. The proof of convergence requires a lemma (Lemma 4.2) which is of some difficulty, because the system (1.2) is infinite dimensional. Its proof is the content of Section 5.

The results of our paper give a rather complete description of this infinite system (1.2) of ordinary differential equations. Similar problems have recently been studied in other contexts, such as coagulation–fragmentation equations [2, 18], although for systems of equations of rather different structure. It is possible that our methods could be useful in other contexts as well.

We conclude the introduction by outlining the results that we need from [1]. First, note that the system (1.2) can be written in a more compact way as

$$
\begin{aligned}
p'_i = {}& -(\lambda_i + \mu_i)p_i + \lambda_{i-1}p_{i-1} + \mu_{i+1}p_{i+1} \\
& + \rho\gamma\left(\sum_{j=0}^{\infty} jp_j\right)(p_{i-1} - p_i) + \nu\left(\delta_{i0}\sum_{j=0}^{\infty} p_j - p_i\right),
\end{aligned}
\tag{1.3}
$$

where

$$p_{-1}(t) := 0 \quad \text{for all } t, \quad \mu_0 := 0 \quad \text{and} \quad \lambda_0 := 0;$$

$$\lambda_i := b_i i \quad \text{and} \quad \mu_i := (d_i + \gamma)i \quad \text{for all } i \geq 1.$$

It is proved in [1] that (1.3) is a well-posed problem in the space $m^1$ defined by

$$m^1 = \left\{ x = (x_0, x_1, \ldots)^T \in \ell^1, \ \sum_j j|x_j| < \infty \right\},$$

equipped with the norm

$$\|\underline{x}\|_m = |x_0| + \sum_{i=0}^{\infty} i|x_i|.$$

More precisely, if **Q** is the infinite matrix

$$(\mathbf{Q})_{ij} = q_{ij} = \begin{cases} b_i i, & \text{if } i+1 = j > 0, \\ -((b_i + (d_i + \gamma))i + \nu), & \text{if } i = j, \\ (d_i + \gamma)i, & \text{if } i - 1 = j > 0, \\ \nu(1 - \delta_{i0}) + (d_1 + \gamma)\delta_{i1}, & \text{if } j = 0, \\ 0, & \text{otherwise,} \end{cases}$$

and $q_j = -q_{jj}$, we define the operator $A$ by

(1.4)
$$\mathcal{D}(A) = \left\{ u \in m^1 : \sum_k q_k |u_k| < \infty \right\};$$

$$Au = u\mathbf{Q}; \qquad (Au)_i = \sum_k q_{ki} u_k.$$

Then it turns out that the closure $\bar{A}$ of $A$ is the generator of a $C^0$-semigroup on $m^1$ (see also [24]).

We then define the map $F : m^1 \to m^1$ by

(1.5)
$$F(\underline{p}) = \rho\gamma \left( \sum_{j=0}^{\infty} j p_j \right) (T_{-1}(\underline{p}) - I(\underline{p})),$$

where $(T_{-1}(\underline{p}))_i := p_{i-1}$ and $I$ is the identity. $F$ is Lipschitz and, in this notation, the system (1.3) can be written as

(1.6)
$$\underline{p}' = A\underline{p} + F(\underline{p});$$
$$\underline{p}(0) = \underline{p}^0.$$

The following theorem is proved in [1].

THEOREM A. *For every $\underline{p}^0 \geq 0 \in \mathcal{D}(\bar{A})$ and any $T > 0$, there exists a unique $\underline{p}(t) \geq 0 \in C([0,T]; \mathcal{D}(\bar{A})) \cap C^1([0,T]; m^1)$ satisfying (1.6). Clearly, $\underline{p}(t)$ will also satisfy (1.3) componentwise.*

*Moreover, if $\underline{p}^0 \in C = \{p \in m^1 : p \geq 0, \sum_{j=0}^{\infty} p_j = 1\}$, then $p(t) \in C$ for all $t \geq 0$.*

Since the $p_i(t)$ represent the frequencies of sites with $i$ individuals, the condition $p(t) \in C$ is quite natural, and most of the following results relate only to that case.

**2. Immigration, birth, death and catastrophe processes.** The analysis of the differential equations system (1.2) is accomplished indirectly, using the properties of a number of associated birth and death processes. We make several comparisons based on couplings of such processes, which exploit the fact that birth and death processes cannot cross without meeting. A good general reference is [20]; in particular, see pages 3 and 4. We begin with a simple lemma.

LEMMA 2.1.  *Fix a positive integer $J$, and let $V = (V_t,\ t \geq 0)$ be the birth and death process on the integers $j \geq J$ with transition rates*

(2.1)
$$j \to j+1 \quad \text{at rate} \quad j\phi,\ j \geq J;$$
$$j \to j-1 \quad \text{at rate} \quad j\mu,\ j \geq J+1,$$

*for some $\phi, \mu > 0$. Then, if $\mathbb{E}^{(m)}$ denotes expectation conditional on $V_0 = m$, for any $j' \geq J$, we have the following:*

1. *If $\phi < \mu$, then $\mathbb{E}^{(j')}\{V_t^2\} \leq \{j'\mu/(\mu-\phi)\}^2$.*
2. *If $\phi > \mu$, then $\mathbb{E}^{(j')}\{V_t^2\} \leq \{2j'\phi e^{(\phi-\mu)t}/(\phi-\mu)\}^2$.*
3. *If $\phi = \mu$, then, for any $\varepsilon > 0$, $\mathbb{E}^{(j')}\{V_t^2\} \leq \{2j'(\phi+\varepsilon)e^{\varepsilon t}/\varepsilon\}^2$.*

PROOF. It is enough to conduct the proof for $j' = J$: for $j' > J$, the $V$-process is stochastically smaller than a $V$-process defined with $J$ replaced by $j'$.

Suppose first that $\phi < \mu$, in which case $V$ is positive recurrent. Observe that a monotone coupling of two realizations of $V$-processes, one with initial state $J$ and the other starting with its equilibrium distribution $\pi$, shows that $\mathbb{E}^{(J)}(V_t^2) \leq \mathbb{E}^\pi(V_0^2)$ for all $t$. Now $\pi$ satisfies the detailed balance equation

$$j\phi\pi(j) = (j+1)\mu\pi(j+1), \qquad j \geq J;$$

hence, $j\pi(j) \leq J(\phi/\mu)^{j-J}$ for all $j \geq J$, from which it follows that

$$\mathbb{E}^{(J)}(V_t^2) \leq \mathbb{E}^\pi(V_0^2) \leq J\sum_{j \geq J} j(\phi/\mu)^{j-J}$$
$$= J\phi\mu(\mu-\phi)^{-2} + J^2\mu(\mu-\phi)^{-1};$$

this proves part 1.

If $\phi > \mu$, we have

(2.2) $$\mathbb{E}V_t = \mathbb{E}(V_t I[\tau_1 \leq t]) + \mathbb{E}(V_t I[\tau_1 > t]),$$

where

$$\tau_1 = \inf\left\{t > 0;\ V_t = J,\ \max_{0 \leq s \leq t} V_s \geq J+1\right\} \leq \infty.$$

Note that $\mathbb{P}^{(J)}[\tau_1 < \infty] = \mu/\phi$ and that $\mathbb{E}^{(J)}(V_t^2 I[\tau_1 > t]) \leq \mathbb{E}^{(J)}(\tilde{V}_t^2)$, where $\tilde{V}$ is a birth and death process on $\mathbb{Z}_+$ with rates as in (2.1), but now for all

$j \geq 0$; this latter bound implies that

$$\mathbb{E}^{(J)}(V_t^2 I[\tau_1 > t]) \leq \{\mathbb{E}^{(J)} \tilde{V}_t\}^2 + \mathrm{Var}^{(J)} \tilde{V}_t$$
(2.3)
$$\leq \{Je^{(\phi-\mu)t}\}^2 + J\,\mathrm{Var}^{(1)} \tilde{V}_t$$
$$\leq J^2 e^{2(\phi-\mu)t} + Je^{2(\phi-\mu)t}\{2(\phi+\mu)/(\phi-\mu) - 1\}.$$

Also, again by a monotone coupling of two $V$-processes,

$$\mathbb{E}^{(J)}(V_t^2 I[\tau_1 \leq t]) \leq \mathbb{P}^{(J)}[\tau_1 \leq t]\mathbb{E}^{(J)}(V_t^2) = (\mu/\phi)\mathbb{E}^{(J)}(V_t^2),$$

and hence, from (2.2) and (2.3),

$$\mathbb{E}^{(J)}(V_t^2) \leq \{\phi/(\phi-\mu)\} J^2 e^{2(\phi-\mu)t}\{2(\phi+\mu)/(\phi-\mu)\},$$

proving part 2, and also, once more by stochastic comparison, part 3. □

Now let $X := (X_t, t \geq 0)$ be an immigration, birth and death process with per capita birth and death rates $\beta_j$ and $\delta_j$, respectively, $j \geq 1$, and with immigration rate $\lambda$. Suppose that the function $n\beta_n$ is concave and increasing in $n \geq 0$, and that $n\delta_n$ is convex and increasing. Then it follows, in particular, that $\beta_n$ is decreasing and $\delta_n$ is increasing in $n \geq 1$; we define

(2.4) $$c := \lim_{n \to \infty} \beta_n - \lim_{n \to \infty} \delta_n.$$

THEOREM 2.2. *Let $X$ and $c$ be as above. Then:*
1. *There exist constants $C_1$ and $C_1(\varepsilon)$, $\varepsilon > 0$, such that*

$$\mathbb{E}^{(j)}(X_t^2) \leq \begin{cases} C_1(1+j^2), & \text{if } c < 0, \\ C_1(\varepsilon)e^{2(c+\varepsilon)t}(1+j^2), & \text{for any } \varepsilon > 0, \text{ if } c \geq 0. \end{cases}$$

2. *There exist constants $C_2$ and $C_2(\varepsilon)$, $\varepsilon > 0$, such that, for all $m \geq 0$,*

$$0 \leq \mathbb{E}^{(m+1)} X_t - \mathbb{E}^{(m)} X_t \leq \begin{cases} C_2, & \text{if } c < 0, \\ C_2(\varepsilon)e^{(c+\varepsilon)t}, & \text{for any } \varepsilon > 0, \text{ if } c \geq 0. \end{cases}$$

3. *In either case, for all $m \geq 0$,*

$$\mathbb{E}^{(m+1)} X_t - \mathbb{E}^{(m)} X_t > \mathbb{E}^{(m+2)} X_t - \mathbb{E}^{(m+1)} X_t.$$

PROOF. Let $\beta'_j := \beta_j + j^{-1}\lambda$ for $j \geq 1$. Then note that, for any positive integer $J$, a simple monotone coupling of two birth and death processes shows that, if $X_0 \leq J$, then $X$ is stochastically smaller than a birth and death process $V$ as in Lemma 2.1, having $\phi = \beta'_J$ and $\mu = \delta_J$ and starting with $V_0 = J$, since $V_0 \geq X_0$ and the sequences $\beta'_j$ and $\delta_j$ are nonincreasing and nondecreasing in $j$, respectively. If $c < 0$, choose $J$ so that $\beta'_J < \delta_J$, and use Lemma 2.1 part 1 to give

$$
(2.5) \quad \begin{aligned} \mathbb{E}^{(j)}(X_t^2) &\leq \mathbb{E}^{(J)}(V_t^2) \leq \{J\delta_J/(\delta_J - \beta'_J)\}^2, \quad j \leq J, \\ \mathbb{E}^{(j)}(X_t^2) &\leq \mathbb{E}^{(j)}(V_t^2) \leq \{j\delta_J/(\delta_J - \beta'_J)\}^2, \quad j > J. \end{aligned}
$$

If $c \geq 0$, choose $J$ so that

$$\delta_J < \beta'_J < \delta_J + c + \varepsilon,$$

if this can be done, and use Lemma 2.1 part 2 as above to give

$$(2.6) \quad \mathbb{E}^{(j)}(X_t^2) \leq \mathbb{E}^{(j)}(V_t^2) \leq \{2\max\{J,j\}\beta'_J e^{(c+\varepsilon)t}/(\beta'_J - \delta_J)\}^2.$$

The only remaining case occurs when $\lambda = 0$ and the sequences $\beta_j$ and $\delta_j$ are both constant for all $j \geq J$ for some $J$, in which case Lemma 2.1 part 2 or 3 can be applied directly. Combining this observation with (2.5) and (2.6), part 1 is proved.

We now turn to part 2. Let $(Y,W) := ((Y_t, W_t), t \geq 0)$ be a two-dimensional pure jump Markov process with transition rates given by

$$
(2.7) \quad \begin{aligned} (i,j) &\to (i+1,j) \quad \text{at rate} \quad i\beta_i + \lambda, \\ (i,j) &\to (i-1,j) \quad \text{at rate} \quad i\delta_i, \\ (i,j) &\to (i,j+1) \quad \text{at rate} \quad (i+j)\beta_{i+j} - i\beta_i, \\ (i,j) &\to (i,j-1) \quad \text{at rate} \quad (i+j)\delta_{i+j} - i\delta_i, \end{aligned}
$$

for all $i, j \geq 0$. All the transition rates are nonnegative, because both $n\beta_n$ and $n\delta_n$ are increasing. Then the processes $Y$ and $Y + W$ are also Markovian, both having the same generator as the immigration, birth and death process $X$. Thus, we can couple realizations $X^1$ and $X^2$ of $X$ with $X_0^1 = m$ and $X_0^2 = m + 1$ by realizing $(Y,W)$ with $Y_0 = m$ and $W_0 = 1$, and setting $X^1 = Y$ and $X^2 = Y + W$. Then it is immediate that $X_t^2 - X_t^1 = W_t \geq 0$ for all $t$; the next step is to bound $\mathbb{E}W_t$.

However, just as before, a simple monotone coupling shows that $W$ is stochastically smaller than a birth and death process $V$ as in Lemma 2.1, having $\phi = \beta'_J$ and $\mu = \delta_J$ and starting with $V_0 = J$, since $V_0 \geq W_0$ and, for any $i \geq 0$,

$$(i+j)\beta_{i+j} - i\beta_i \leq j\beta_j \leq j\beta'_J, \qquad j \geq J,$$

and

$$(i+j)\delta_{i+j} - i\delta_i \geq j\delta_j \geq j\delta_J, \qquad j \geq J+1,$$

by the concavity of $n\beta_n$ and the convexity of $n\delta_n$. Thus, in particular, $\mathbb{E}W_t \leq \mathbb{E}^{(J)}V_t$, and the bounds on $\mathbb{E}^{(J)}(V_t^2)$ obtained in part 1 can be invoked, completing the proof of part 2.

A METAPOPULATION MODEL 9

For part 3, we extend $(Y, W)$ to a four-dimensional pure jump Markov process $((Y_t, W_t, U_t, V_t), \ t \geq 0)$ with transition rates

$$\mathbf{n} \to \mathbf{n} + \boldsymbol{\epsilon}^{(1)} \quad \text{at rate} \quad i\beta_i + \lambda,$$
$$\mathbf{n} \to \mathbf{n} - \boldsymbol{\epsilon}^{(1)} \quad \text{at rate} \quad i\delta_i,$$
$$\mathbf{n} \to \mathbf{n} + \boldsymbol{\epsilon}^{(2)} \quad \text{at rate} \quad (i+j)\beta_{i+j} - i\beta_i,$$
$$\mathbf{n} \to \mathbf{n} - \boldsymbol{\epsilon}^{(2)} \quad \text{at rate} \quad (i+j)\delta_{i+j} - i\delta_i,$$
$$\mathbf{n} \to \mathbf{n} + \boldsymbol{\epsilon}^{(3)} \quad \text{at rate} \quad (i+k)\beta_{i+k} - i\beta_i,$$
$$\mathbf{n} \to \mathbf{n} - \boldsymbol{\epsilon}^{(3)} \quad \text{at rate} \quad (i+k)\delta_{i+k} - i\delta_i,$$
$$\mathbf{n} \to \mathbf{n} + \boldsymbol{\epsilon}^{(4)} \quad \text{at rate} \quad (i+j+l)\beta_{i+j+l} - (i+j)\beta_{i+j},$$
$$\mathbf{n} \to \mathbf{n} - \boldsymbol{\epsilon}^{(4)} \quad \text{at rate} \quad (i+j+l)\delta_{i+j+l} - (i+j)\delta_{i+j},$$

when $\mathbf{n} = (i, j, k, l)$ is such that $k \neq l$, the last four transitions being replaced by

$$\mathbf{n} \to \mathbf{n} + \boldsymbol{\epsilon}^{(3)} + \boldsymbol{\epsilon}^{(4)} \quad \text{at rate} \quad (i+j+k)\beta_{i+j+k} - (i+j)\beta_{i+j},$$
$$\mathbf{n} \to \mathbf{n} - \boldsymbol{\epsilon}^{(3)} - \boldsymbol{\epsilon}^{(4)} \quad \text{at rate} \quad (i+k)\delta_{i+k} - i\delta_i,$$

$$\mathbf{n} \to \mathbf{n} + \boldsymbol{\epsilon}^{(3)} \quad \text{at rate} \quad (i+k)\beta_{i+k} - i\beta_i - (i+j+k)\beta_{i+j} + (i+j)\beta_{i+j},$$
$$\mathbf{n} \to \mathbf{n} - \boldsymbol{\epsilon}^{(4)} \quad \text{at rate} \quad (i+j+k)\delta_{i+j+k} - (i+j)\delta_{i+j} - (i+k)\delta_{i+k} + i\delta_i,$$

when $\mathbf{n} = (i, j, k, k)$, all transition rates being nonnegative because of the assumptions on $nb_n$ and $nd_n$; here, $\boldsymbol{\epsilon}^{(m)}$ denotes the $m$th coordinate vector. The four processes $Y$, $Y+W$, $Y+U$ and $Y+W+V$ are Markov, and each has the same generator as the immigration, birth and death process $X$. Thus, realizations $X^1$, $X^2$, $X^3$ and $X^4$ of $X$ with $X_0^1 = m$, $X_0^2 = X_0^3 = m+1$ and $X_0^4 = m+2$ can be obtained from $(Y, W, U, V)$ by setting $X^1 = Y$, $X^2 = Y + U$, $X^3 = Y + W$ and $X^4 = Y + W + V$ and taking $Y_0 = m$, $W_0 = U_0 = V_0 = 1$. Thus, $\mathbb{E}^{(m+1)} X_t - \mathbb{E}^m X_t = \mathbb{E} U_t$ and $\mathbb{E}^{(m+2)} X_t - \mathbb{E}^{(m+1)} X_t = \mathbb{E} V_t$. Initially, $U_0 = V_0 = 1$. Thereafter, both $U$ and $V$ make only unit jumps, and at any time at which $U$ and $V$ are equal, either they can jump together, or $U$ can increase by 1 or $V$ can decrease by 1. Thus, $U$ is always greater than or equal to $V$, and, for each $t > 0$, $U_t > V_t$ with positive probability. Hence, for all $t > 0$,

$$\mathbb{E}^{(m+1)} X_t - \mathbb{E}^m X_t = \mathbb{E} U_t > \mathbb{E} V_t = \mathbb{E}^{(m+2)} X_t - \mathbb{E}^{(m+1)} X_t,$$

proving part 3. $\square$

The theorem above is used in the study of our main object of interest, a family of immigration, birth, death and catastrophe processes $Z^{(s)}$, indexed by an immigration parameter $s$. The pure jump Markov process $Z^{(s)}$ has transition rates

$$
\begin{aligned}
j \to j+1 & \quad \text{at rate} \quad q_{j,j+1} := jb_j + \rho\gamma s, \\
j \to j-1 & \quad \text{at rate} \quad q_{j,j-1} := j(d_j + \gamma), \\
j \to 0 & \quad \text{at rate} \quad q_{j,0} := \nu,
\end{aligned}
\tag{2.8}
$$

and $nb_n$ is assumed to be increasing and concave, $nd_n$ to be increasing and convex. The process $Z^{(s)}$ starting with any initial distribution $\psi$ can be constructed as follows from a sequence of *independent* realizations $X^{(0)}, X^{(1)}, \ldots$ of an $X$-process with parameters $\beta_i = b_i$, $\lambda = \rho\gamma s$ and $\delta_i = d_i + \gamma$, and with $X_0^{(0)} \sim \psi$ and $X_0^{(n)} = 0$, $n \geq 1$. Let the times $(T_n, n \geq 1)$ of the catastrophes be the partial sums of independent negative exponentially distributed random variables $(E_n, n \geq 1)$ with mean $1/\nu$, which are also independent of $(X^{(n)}, n \geq 0)$. Set $N(t) := \min\{n \geq 0 : T_n \leq t\}$, where $T_0 := 0$; then define

$$Z_t^{(s)} := X^{(N(t))}(t - T_{N(t)}). \tag{2.9}$$

A pair of $Z^{(s)}$-processes $Z^{(s,1)}$ and $Z^{(s,2)}$ with different initial states $k > l$ can then always be coupled by using the same sequence of $X$-processes $(X^{(n)}, n \geq 1)$ and taking $X^{(0,1)} = Y + W$, $X^{(0,2)} = Y$, where $(Y, W)$ is as in (2.7) and has $Y_0 = l$ and $W_0 = k - l$. With this construction, it is clear that $Z_t^{(s,1)} \geq Z_t^{(s,2)}$ for all $t$, that $\mathbb{P}[Z_t^{(s,1)} > Z_t^{(s,2)}] \leq e^{-\nu t}$ and that

$$0 \leq \mathbb{E}(Z_t^{(s,1)} - Z_t^{(s,2)}) = e^{-\nu t}\mathbb{E}W_t = e^{-\nu t}(\mathbb{E}^{(k)}X_t - \mathbb{E}^{(l)}X_t). \tag{2.10}$$

Defining

$$c := \lim_{n \to \infty} b_n - \lim_{n \to \infty} d_n - \gamma, \tag{2.11}$$

and assuming that $c < \nu$, it thus follows from Theorem 2.2 part 2 that, for $k > l$,

$$
\begin{aligned}
0 &\leq \mathbb{E}^{(k)}Z_t^{(s)} - \mathbb{E}^{(l)}Z_t^{(s)} \\
&\leq \begin{cases} (k-l)C_2 e^{-\nu t}, & \text{if } c < 0, \\ (k-l)C_2(\tfrac{1}{2}(\nu - c))\exp\{-\tfrac{1}{2}(\nu - c)t\}, & \text{if } c \geq 0. \end{cases}
\end{aligned}
\tag{2.12}
$$

Thus, if $f : \mathbb{Z}_+ \to \mathbb{R}$ is Lipschitz with constant $K(f)$, then

$$|\mathbb{E}^{(k)}f(Z_t^{(s)}) - \mathbb{E}^{(l)}f(Z_t^{(s)})| \leq C(k-l)K(f)e^{-\alpha t} \tag{2.13}$$

for some $C, \alpha > 0$. Furthermore, from (2.10) and Theorem 2.2 part 3, we have

$$\mathbb{E}^{(m+1)}Z_t^{(s)} - \mathbb{E}^{(m)}Z_t^{(s)} > \mathbb{E}^{(m+2)}Z_t^{(s)} - \mathbb{E}^{(m+1)}Z_t^{(s)}, \tag{2.14}$$

for all $m, t \geq 0$.

THEOREM 2.3. *Let $Z^{(s)}$ be as defined in (2.8), with $nb_n$ increasing and concave, $nd_n$ increasing and convex. Suppose that $c < \nu$, where $c$ is as defined in (2.11). Then, for $s > 0$, $Z^{(s)}$ is positive recurrent, and its equilibrium distribution $\pi^{(s)}$ has finite mean equal to $\lim_{t \to \infty} \mathbb{E}^{(0)} Z_t^{(s)}$; furthermore, for any $0 \leq \delta \leq 1$ for which $c(1+\delta) < \nu$, we can find $K_1(\delta) < \infty$ such that*

$$\mathbb{E}^{(j)}\{(Z_t^{(s)})^{(1+\delta)}\} \leq K_1(\delta)\{1 + j^{(1+\delta)}\}, \tag{2.15}$$

*for all $t \geq 0$ and $j \geq 0$.*

*If $s = 0$, the state $0$ is absorbing for $Z^{(s)}$, and the only stationary distribution is $\pi^{(0)} = \Delta_{\{0\}}$, giving probability one to $0$.*

PROOF. The case $s = 0$ is immediate, so we now suppose that $s > 0$.

If $\nu = 0$ and $c < 0$, the detailed balance equations

$$(jb_j + \rho\gamma s)\pi_j = (j+1)(d_{j+1} + \gamma)\pi_{j+1}, \qquad j \geq 0, \tag{2.16}$$

are satisfied with

$$\pi_{j+1} \leq \frac{b_j + \rho\gamma s j^{-1}}{d_{j+1} + \gamma}\pi_j \leq (1-\varepsilon)\pi_j,$$

for some $\varepsilon > 0$ and for all $j$ large enough, because $c < 0$. Hence, (2.16) have a nonnegative solution with geometrically decreasing tail, and the conclusion of the theorem follows.

If $\nu > 0$, positive recurrence is immediate. Construct $Z^{(s)}$ with $Z_0^{(s)} = 0$ from a sequence of $X$-processes as above. Then, if $m(t) := m^{(s)}(t) := \mathbb{E} Z_t^{(s)}$, we have the renewal equation

$$m(t) = e^{-\nu t}\mathbb{E}\{X_t | X_0 = 0\} + \int_0^t \nu e^{-\nu u} m(t-u)\, du.$$

Now, by Lemma 2.1 part 1,

$$\mathbb{E}^{(0)} X_t \leq \begin{cases} C_1, & \text{if } c < 0, \\ C_2((\nu-c)/2)\exp\{\tfrac{1}{2}(\nu+c)t\}, & \text{if } c \geq 0, \end{cases}$$

for suitable constants $C_1$ and $C_2$. Furthermore, a monotone coupling of two $X$-processes with different initial conditions shows that $\mathbb{E}^{(0)} X_t$ increases with $t$. Hence, the key renewal theorem ([11], page 363) can be applied to conclude that

$$m^{(s)}(\infty) := \lim_{t \to \infty} \mathbb{E} Z_t^{(s)} = \nu \int_0^\infty e^{-\nu t}\mathbb{E}\{X_t | X_0 = 0\}\, dt \tag{2.17}$$

exists and is finite. But now, because $Z^{(s)}$ is nonnegative and positive recurrent, it follows from (2.17) that $\pi^{(s)}$ has finite mean, satisfying

$$\pi^{(s)}(e) = \mathbb{E}^{\pi^{(s)}}(Z_0^{(s)}) \leq m^{(s)}(\infty),$$

where $e(j) := j$ for all $j \geq 0$ and $\pi^{(s)}(f) := \sum_{k \geq 0} \pi_k^{(s)} f(k)$.

Finally, for any $0 < \delta \leq 1$ for which $c(1 + \delta) < \nu$, a similar renewal argument can be employed, again appealing to Lemma 2.1 part 1, to show that

$$m_\delta^{(s)}(t) := \mathbb{E}^{(0)}\{(Z_t^{(s)})^{(1+\delta)}\}$$

is uniformly bounded for all $t$; hence, the sequence of random variables $Z_t^{(s)}$ is uniformly integrable, and thus, in fact,

(2.18) $$\pi^{(s)}(e) = m^{(s)}(\infty),$$

proving the first two claims of the theorem. Noting also that, for any $\varepsilon > 0$,

(2.19) $$\begin{aligned}\mathbb{E}^{(j)}&\{(Z_t^{(s)})^{(1+\delta)}\}\\ &= e^{-\nu t}\mathbb{E}^{(j)}\{X_t^{(1+\delta)}\} + \int_0^t \nu e^{-\nu u} m_\delta^{(s)}(t-u)\, du\\ &\leq m_\delta^{(s)}(t) + (1 + j^{(1+\delta)})C_1^{(1+\delta)/2}(\varepsilon)\exp\{(1+\delta)(c+\varepsilon)t - \nu t\},\end{aligned}$$

from Theorem 2.2 part 1, the remaining claim is also proved. □

With these preparations, we can now prove the main result of the section. The assumptions of Theorem 2.3 are still in force.

THEOREM 2.4. *Let $\pi^{(s)}$ denote the equilibrium distribution of the process $Z^{(s)}$; then $\pi^{(s)}(f)$ is continuous in $s$ for any Lipschitz function $f$. Furthermore, if $e(j) := j$ for all $j \geq 0$, then $\pi^{(s)}(e)$ is an increasing, strictly concave function of $s$.*

PROOF. Let $\mathcal{A}^{(s)}$ be the generator of the process $Z^{(s)}$, so that

(2.20) $$\begin{aligned}(\mathcal{A}^{(s)}h)(j) &:= \sum_{l \neq j} q_{jl}\{h(l) - h(j)\}\\ &= (jb_j + \rho\gamma s)\{h(j+1) - h(j)\}\\ &\quad + j(d_j + \gamma)\{h(j-1) - h(j)\} + \nu\{h(0) - h(j)\},\end{aligned}$$

and, for any Lipschitz function $f$ with constant $K(f)$, let $\theta^{(s)}(f)$ be defined by

(2.21) $$\theta^{(s)}(f)(j) := -\int_0^\infty \{\mathbb{E}^{(j)}f(Z_t^{(s)}) - \pi^{(s)}(f)\}\, dt.$$

We begin by showing that $\theta^{(s)}(f)$ is a solution $h$ to the equation

(2.22) $$(\mathcal{A}^{(s)}h)(j) = f(j) - \pi^{(s)}(f), \qquad j \geq 0.$$

First, realizing $\pi^{(s)}(f) = \mathbb{E}^{\pi^{(s)}} f(Z_t^{(s)})$, it follows from (2.13) that

$$|\mathbb{E}^{(j)} f(Z_t^{(s)}) - \pi^{(s)}(f)|$$

(2.23)
$$= \left| \sum_{k \geq 0} \pi_k^{(s)} \{ \mathbb{E}^{(j)} f(Z_t^{(s)}) - \mathbb{E}^{(k)} f(Z_t^{(s)}) \} \right|$$

$$\leq CK(f) e^{-\alpha t} \sum_{k \geq 0} \pi_k^{(s)} |k - j| \leq CK(f) e^{-\alpha t} \{ m^{(s)}(\infty) + j \}.$$

Hence, $\theta^{(s)}(f)$ given in (2.21) is well defined. Now set

$$\theta_T^{(s)}(f)(j) := - \int_0^T \{ \mathbb{E}^{(j)} f(Z_t^{(s)}) - \pi^{(s)}(f) \} \, dt,$$

noting that $\lim_{T \to \infty} \theta_T^{(s)}(f)(j) = \theta^{(s)}(f)(j)$ by (2.23). Conditioning on the first jump gives

(2.24)
$$\theta_T^{(s)}(f)(j) = -\mathbb{E}^{(j)} \left\{ \int_0^T \{ f(Z_t^{(s)}) - \pi^{(s)}(f) \} \, dt \right\}$$

$$= -e^{-q_j T} T \{ f(j) - \pi^{(s)}(f) \}$$

$$- \int_0^T e^{-q_j u} \sum_{l \neq j} q_{jl} \left\{ \mathbb{E}^{(l)} \left( \int_0^{T-u} \{ f(Z_t^{(s)}) - \pi^{(s)}(f) \} \, dt \right) \right.$$

$$\left. + u \{ f(j) - \pi^{(s)}(f) \} \right\} du$$

$$= -q_j^{-1}(1 - e^{-q_j T}) \{ f(j) - \pi^{(s)}(f) \}$$

$$+ \int_0^T e^{-q_j u} \sum_{l \neq j} q_{jl} \theta_{T-u}^{(s)}(f)(l) \, du.$$

Now, by (2.23),

$$|\theta_{T-u}^{(s)}(f)(l) \mathbb{1}_{\{u \leq T\}}| \leq CK(f) \alpha^{-1} (m^{(s)}(\infty) + l)$$

for all $T$ and $u$, and $q_{jl} > 0$ only for $l = 0, j-1, j+1$. Thus, letting $T \to \infty$ in (2.24) and using dominated convergence, it follows that

$$\theta^{(s)}(f)(j) = -q_j^{-1} \{ f(j) - \pi^{(s)}(f) \} + \sum_{l \neq j} q_{jl} \int_0^\infty e^{-q_j u} \theta^{(s)}(f)(l) \, du,$$

or

$$f(j) - \pi^{(s)}(f) = \sum_{l \neq j} q_{jl} \{ \theta^{(s)}(f)(l) - \theta^{(s)}(f)(j) \}.$$

Thus, $\theta^{(s)}(f)$ solves (2.22).

Furthermore, again using (2.23),

$$(2.25) \quad \theta^{(s)}(f)(j+1) - \theta^{(s)}(f)(j) = \int_0^\infty \{\mathbb{E}^{(j)} f(Z_t^{(s)}) - \mathbb{E}^{(j+1)} f(Z_t^{(s)})\} dt$$

and (2.13) immediately gives

$$(2.26) \quad |\theta^{(s)}(f)(j+1) - \theta^{(s)}(f)(j)| \leq CK(f)/\alpha;$$

thus, $\Delta \theta^{(s)}(f)$ is bounded and, hence, also Lipschitz, with constant

$$(2.27) \quad K(\Delta \theta^{(s)}(f)) \leq 2CK(f)/\alpha.$$

Now, by Dynkin's formula ([14], Theorem 2), it follows that $\pi^{(s)}(\mathcal{A}^{(s)} h) = 0$ for all $s$, for any Lipschitz function $h$. In particular, for any $t > -s$, using (2.20),

$$
\begin{aligned}
0 &= \pi^{(s+t)}(\mathcal{A}^{(s+t)} \theta^{(s)}(f)) = \mathbb{E}^{\pi^{(s+t)}}(\mathcal{A}^{(s+t)} \theta^{(s)}(f))(Z_0^{(s)}) \\
(2.28) \quad &= \mathbb{E}^{\pi^{(s+t)}} \{(\mathcal{A}^{(s)} \theta^{(s)}(f))(Z_0^{(s)}) + \rho \gamma t (\Delta \theta^{(s)}(f))(Z_0^{(s)})\} \\
&= \mathbb{E}^{\pi^{(s+t)}} \{f(Z_0^{(s)}) - \pi^{(s)}(f) + \rho \gamma t (\Delta \theta^{(s)}(f))(Z_0^{(s)})\}.
\end{aligned}
$$

Thus, from (2.28) and (2.26), it follows that

$$(2.29) \quad |\pi^{(s+t)}(f) - \pi^{(s)}(f)| \leq \rho \gamma t \|\Delta \theta^{(s)}(f)\| \leq \rho \gamma |t| CK(f)/\alpha$$

for any Lipschitz function $f$, so that $\pi^{(s)}(f)$ is continuous in $s$, proving the first part of the theorem.

It then also follows that

$$
\begin{aligned}
|\pi^{(s+t)}(f) &- \pi^{(s)}(f) + \rho \gamma t \pi^{(s)}(\Delta \theta^{(s)}(f))| \\
&\leq \rho \gamma |t| |\pi^{(s+t)}(\Delta \theta^{(s)}(f)) - \pi^{(s)}(\Delta \theta^{(s)}(f))|
\end{aligned}
$$

and, hence, that

$$(2.30) \quad \frac{d}{ds} \pi^{(s)}(f) = -\rho \gamma \pi^{(s)}(\Delta \theta^{(s)}(f)).$$

Taking $f = e$, this last can be re-expressed using (2.21) as

$$(2.31) \quad \frac{d}{ds} \pi^{(s)}(e) = \rho \gamma \mathbb{E}^{\pi^{(s)}} \int_0^\infty \{g(Z_0^{(s)} + 1, t) - g(Z_0^{(s)}, t)\} dt,$$

where $g(j, t) := \mathbb{E}^{(j)} Z_t^{(s)}$. Hence, from Theorem 2.2 part 2, it follows that $\pi^{(s)}(e)$ is increasing in $s$, proving the next part of the theorem.

Now, from (2.28) with $u$ and $2u$ for $t$,

$$\pi^{(s+u)}(f) - \pi^{(s)}(f) = -\rho \gamma u \pi^{(s+u)}(\Delta \theta^{(s)}(f))$$

and
$$\pi^{(s+2u)}(f) - \pi^{(s)}(f) = -2\rho\gamma u \pi^{(s+2u)}(\Delta\theta^{(s)}(f))$$

giving, again from (2.28),

$$\pi^{(s+2u)}(f) - 2\pi^{(s+u)}(f) + \pi^{(s)}(f)$$
$$= -2\rho\gamma u\{\pi^{(s+2u)}(\Delta\theta^{(s)}(f)) - \pi^{(s+u)}(\Delta\theta^{(s)}(f))\}$$
$$= -2\rho\gamma u\{\pi^{(s)}(\Delta\theta^{(s)}(f)) - 2\rho\gamma u \pi^{(s+2u)}(\Delta\theta^{(s)}(\Delta\theta^{(s)}(f)))$$
$$\quad - \pi^{(s)}(\Delta\theta^{(s)}(f)) + \rho\gamma u \pi^{(s+u)}(\Delta\theta^{(s)}(\Delta\theta^{(s)}(f)))\}$$
$$= 2(\rho\gamma u)^2 \pi^{(s)}(\Delta\theta^{(s)}(\Delta\theta^{(s)}(f))) + \eta,$$

where

$$|\eta| \le 10(\rho\gamma|u|)^3 CK(\Delta\theta^{(s)}(\Delta\theta^{(s)}(f)))/a \le 40(C\rho\gamma|u|/\alpha)^3 K(f),$$

this last by (2.29) and (2.27). Hence, $\pi^{(s)}(f)$ is twice differentiable, and

(2.32) $$\frac{d^2}{ds^2}\pi^{(s)}(f) = 2\rho^2\gamma^2 \pi^{(s)}(\Delta\theta^{(s)}(\Delta\theta^{(s)}(f))).$$

Now, using the formula given in (2.25), it follows that

$$\Delta\theta^{(s)}(\Delta\theta^{(s)}(f))(m)$$
(2.33)
$$= -\int_0^\infty \{\mathbb{E}^{(m+1)}\Delta\theta^{(s)}(f)(Z_t^{(s)}) - \mathbb{E}^{(m)}\Delta\theta^{(s)}(f)(Z_t^{(s)})\}\,dt$$
$$= \int_0^\infty \left\{\mathbb{E}^{(m+1)}\int_0^\infty \{g_f(Z_t^{(s)}+1,w) - g_f(Z_t^{(s)},w)\}\,dw\right.$$
$$\left. - \mathbb{E}^{(m)}\int_0^\infty \{g_f(Z_t^{(s)}+1,w) - g_f(Z_t^{(s)},w)\}\,dw\right\}dt,$$

where $g_f(l,w) := \mathbb{E}^{(l)} f(Z_w^{(s)})$. To evaluate (2.33), realize $Z^{(s,1)}$ with $Z_0^{(s,1)} = j+1$ and $Z^{(s,2)}$ with $Z_0^{(s,2)} = j$ as before, using the Markov process $(Y,W)$ of (2.7), with $Y_0 = j$ and $W_0 = 1$, so that

$$Z_t^{(s,1)} = Z_t^{(s,2)} + W_t I[E_1 > t],$$

where $E_1$ is an independent negative exponential random variable with mean $1/\nu$. Thus,

$$\Delta\theta^{(s)}(\Delta\theta^{(s)}(f))(j)$$
(2.34) $$= \mathbb{E}^{(j)} \int_0^\infty \int_0^\infty e^{-\nu t}\{g_f(Y_t+W_t+1,w) - g_f(Y_t+W_t,w)$$
$$\quad - g_f(Y_t+1,w) + g_f(Y_t,w)\}\,dw\,dt.$$

In order to use (2.32) to investigate the curvature of $\pi^{(s)}(e)$, we take $f = e$ in (2.34). Then, for any $k \geq l$,

$$g_e(k+1,w) - g_e(k,w) - g_e(l+1,w) + g_e(l,w)$$
$$= \{\mathbb{E}^{(k+1)} Z_w^{(s)} - \mathbb{E}^{(k)} Z_w^{(s)}\} - \{\mathbb{E}^{(l+1)} Z_w^{(s)} - \mathbb{E}^{(l)} Z_w^{(s)}\} < 0,$$

from (2.14), for all $w > 0$, so that the integrand is always negative. Hence, from (2.34), it follows that $\Delta\theta^{(s)}(\Delta\theta^{(s)}(e))(j) < 0$ for all $j$ and $s$, and thus, from (2.32),

$$\frac{d^2}{ds^2}\pi^{(s)}(e) < 0, \qquad s \geq 0.$$

This completes the proof of the theorem. □

**3. Equilibria.** We now investigate the equilibrium solutions of (1.3). For the sake of simplicity, we shall assume here and in all that follows that $\rho = 1$. There is no real loss of generality in this, since one could set $d'_i = d_i + \gamma(1-\rho)$ and $\gamma' = \gamma\rho$, and write (1.2) using $d'$ and $\gamma'$ in place of $d$ and $\gamma$. In biological terms, unsuccessful migration is just one cause of death.

If $\pi \in m^1$ is such a solution, and

$$s := \sum_{j=0}^{\infty} j\pi_j,$$

then $\pi$ must solve

$$0 = -[(b_i + d_i + \gamma)i + \nu + \gamma s]\pi_i$$
(3.1)
$$\qquad + [b_{i-1}(i-1) + \gamma s]\pi_{i-1} + [d_{i+1} + \gamma](i+1)\pi_{i+1}, \qquad i \geq 1;$$
$$0 = \nu(1 - \pi_0) + (d_1 + \gamma)\pi_1 - \gamma s\pi_0.$$

Hence, $\pi$ must be the equilibrium distribution of the immigration, birth, death and catastrophe process $Z_t^{(s)}$, which we studied in detail in Section 2. From Theorem 2.3, and using (H1) and (H2), we know that $Z^{(s)}$ has a unique stationary distribution $\pi^{(s)}$, which has finite mean denoted by

(3.2) $$G(s) := \pi^{(s)}(e) = \sum_{j \geq 1} j\pi_j^{(s)}.$$

In order to have an equilibrium solution of (1.3), $s$ must be equal to $\pi^{(s)}(e)$; in other words, we look for a solution to the equation

(3.3) $$s = G(s),$$

a fixed point of the function $G$.

THEOREM 3.1. *Suppose that* (H1) *and* (H2) *are satisfied. If* $G'(0) > 1$, *then there exists a unique positive fixed point* $s^*$ *of* $G$; *if* $G'(0) \leq 1$, *then* $G(s) < s$ *for all* $s > 0$.

REMARK 3.2. Note that $s = 0$ is always a fixed point of $G$; the corresponding equilibrium distribution $\pi^{(0)}$ is the vector $e^0 = (1, 0, 0, \ldots)^T$, which can be interpreted as the extinction equilibrium.

For the proof, we need a technical point.

LEMMA 3.3. *Let*

$$m^2 = \left\{ x \in \ell^1, \sum_j j^2 |x_j| < \infty \right\},$$

*with norm*

$$\|\underline{x}\|_{m^2} = |x_0| + \sum_{i=0}^{\infty} i^2 |x_i|,$$

*and let* $A_2$ *be the part of* $\bar{A}$ *in* $m^2$, *that is,*

$$\mathcal{D}(A_2) = \{ x \in \mathcal{D}(\bar{A}) : \bar{A}x \in m^2 \}; \qquad A_2 x = \bar{A}x.$$

*Then, if* $p(0) \in \mathcal{D}(A_2)$, $p(t)$ *satisfies*

(3.4) $$\sum_{j=1}^{\infty} j^2 d_j p_j(t) < \infty.$$

PROOF. We first note that the restriction of $e^{\bar{A}t}$ to $m^2$ is again a $C^0$-semigroup. This can be established following, with obvious changes, the proofs in [1]. In fact, repeating step by step the proof of Proposition 6.5 of [1], one sees that $A_2 - \omega$ is dissipative, as long as $\omega \geq 3 \max_i b_i$. The density of the range is then established exactly as in Proposition 6.6 of [1].

Moreover, repeating the proofs of Lemmas 6.8 and 6.9 of [1], one sees that the domain of the restricted semigroup is contained in the set

$$\left\{ x \in \ell^1, \sum_j j^2 d_j |x_j| < \infty \right\},$$

and (3.4) follows. □

PROOF OF THEOREM 3.1. First of all, it is proved in Theorem 2.4 that $G$ is an increasing, strictly concave function in $s \geq 0$. We now establish two further properties of $G$:

(3.5) $$G(0) = 0 \quad \text{and} \quad \lim_{s \to \infty} G(s)/s < 1.$$

The first of these follows because, when $s = 0$, the equilibrium distribution is concentrated at 0, so that its mean is 0.

For the limit as $s \to \infty$, we let $m(t) = m^{(s)}(t) = \mathbb{E}(Z_t^{(s)})$, with $Z_t^{(s)}$ as defined in (2.8), noting that $G(s) = \lim_{t \to \infty} m^{(s)}(t)$ as shown in (2.18). Letting $p_j(t) = \mathbb{P}(Z_t^{(s)} = j)$, we can write $m(t) = \sum_j j p_j(t)$.

The forward equations satisfied by $p(t)$ can be written as $p'(t) = A_s p(t)$, where

$$A_s : \mathcal{D}(\bar{A}) \to m^1, \qquad (A_s p)_i = \begin{cases} (\bar{A}p)_i + \gamma s(p_{i-1} - p_i), & i \geq 1, \\ (\bar{A}p)_0 - \gamma s p_0, & i = 0, \end{cases}$$

is a bounded perturbation of the operator $\bar{A}$ defined in (1.4).

Hence, if the initial value $p(0)$ is in $\mathcal{D}(\bar{A})$, then $p(t) = e^{A_s t} p(0)$ is differentiable as a function from $\mathbb{R}$ to $m^1$ and we have

$$m'(t) = \sum_{j=1}^{\infty} j p'_j(t)$$

(3.6)
$$= \sum_{j=1}^{\infty} j \{ ((j-1)b_{j-1} + \gamma s) p_{j-1}(t)$$
$$- ((b_j + (d_j + \gamma))j + \gamma s + \nu) p_j(t)$$
$$+ (j+1)(d_{j-1} + \gamma) p_{j+1}(t) \}.$$

If $p(0) \in \mathcal{D}(A_2)$, the condition (3.4) allows the order of the sums in (3.6) to be interchanged, and, with some manipulations, we obtain

(3.7) $$m'(t) = \sum_j j(b_j - d_j) p_j(t) - (\gamma + \nu) m(t) + \gamma s.$$

Using the concavity of $xb(x)$ and the convexity of $xd(x)$, we obtain

(3.8) $$\sum_j j b_j p_j(t) \leq m(t) b(m(t)) \quad \text{and} \quad \sum_j j d_j p_j(t) \geq m(t) d(m(t)).$$

Hence, from (3.7), it follows that

(3.9) $$m'(t) \leq m(t)[b(m(t)) - d(m(t)) - \gamma - \nu] + \gamma s,$$

so that $m(t) \leq x^{(s)}(t)$, where $x := x^{(s)}(t)$ is the solution of the Cauchy problem

(3.10)
$$x' = x[b(x) - d(x) - \gamma - \nu] + \gamma s,$$
$$x(0) = m(0).$$

Since $\mathcal{D}(A_2)$ is dense in $\mathcal{D}(A)$, it follows that $m(t) \leq x^{(s)}(t)$ for all $p(0) \in \mathcal{D}(A)$.

Set

(3.11) $$a = \nu + d_\infty - b_\infty > 0$$

because of (H2), and choose $\bar{m}$ such that

$$b(\bar{m}) - d(\bar{m}) - \nu = -a/2,$$

if this is possible; otherwise, set $\bar{m} = 0$. In any case, we have

(3.12) $$b(m) - d(m) - \nu \leq -a/2 \qquad \text{for } m \geq \bar{m}.$$

Take $\bar{s}$ such that $\gamma \bar{s} = \bar{m}(\frac{a}{2} + \gamma)$. Then for all $s > \bar{s}$, there exists $\tau^{(s)}$ such that $x(\tau^{(s)}) = \bar{m}$ and $x(t) > \bar{m}$ for $t > \tau^{(s)}$. Then, using (3.12), we have

$$x'(t) \leq \gamma s - x(t)(\tfrac{1}{2}a + \gamma) \qquad \text{for } t > \tau^{(s)}.$$

Hence,

$$x(t) \leq \bar{m} e^{-(a/2+\gamma)(t-\tau^{(s)})} + \gamma s \int_{\tau^{(s)}}^{t} e^{-(a/2+\gamma)(t-\sigma)}\, d\sigma$$

$$= \frac{\gamma s}{\gamma + a/2} - \left(\frac{\gamma s}{\gamma + a/2} - \bar{m}\right) e^{-(a/2+\gamma)(t-\tau^{(s)})} \leq \frac{\gamma s}{\gamma + a/2},$$

so that, using (2.18),

$$G(s) = \lim_{t\to\infty} m(t) \leq \lim_{t\to\infty} x(t) \leq \frac{\gamma s}{\gamma + a/2}$$

and, hence,

$$\lim_{s\to\infty} G(s)/s \leq \frac{\gamma}{\gamma + a/2} < 1,$$

as stated above.

Turning now to the fixed points of $G$, note that $G(0) = 0$ and $G$ is strictly concave; hence, $G(s) = s$ has, at most, one other solution in $s \geq 0$. Since also $\lim_{s\to\infty} G(s)/s < 1$, it follows that there is a unique positive solution of $G(s)/s = 1$ if $G'(0) > 1$; otherwise, if $G'(0) \leq 1$, we have $G(s)/s < 1$ for all $s > 0$. $\square$

The next result shows that assuming condition (H2) to be satisfied is not restrictive, when looking for positive equilibria of (1.2).

PROPOSITION 3.4. *If* (H2) *is violated, there are no nontrivial equilibrium solutions to* (1.2).

PROOF. If $\gamma = 0$, the proposition follows immediately from the $p^0$-equation in (1.2). Otherwise, the process $Z^{(s)}$ is stochastically larger than a process $\widehat{Z}^{(s)}$ which has $\hat{b}_j = b_\infty$ and $\hat{d}_j = d_\infty$ for all $j$, and the same is true if $\hat{b}_j = b_\infty - b_*$ for any $0 \le b_* < b_\infty$. Letting $\hat{m} := \mathbb{E}^1 \widehat{Z}^{(s)}_t$, note that, as for (3.6),

$$(3.13) \quad \hat{m}'(t) = \hat{m}(t)\{b_\infty - b_* - d_\infty - \gamma - \nu\} + \gamma s = -\hat{m}(t)(a' + \gamma) + \gamma s,$$

where $a' = a + b_*$ and $a$ is as in (3.11). Suppose now that $a < 0$, so that (H2) is violated. If $a \le -\gamma$, choose $b_*$ so that $a' = -\frac{1}{2}\gamma$; otherwise, take $b_* = 0$. Then it follows from (3.13) that

$$G(s) = \lim_{t \to \infty} \mathbb{E} Z^{(s)}_t \ge \lim_{t \to \infty} \hat{m}(t) = s\{\gamma/(\gamma + a')\} > s,$$

for all $s > 0$, and there can be no $s > 0$ for which $G(s) = s$.

Finally, if $a = 0$ and $\gamma > 0$, then the condition $c < \nu$ of Theorem 2.3 is satisfied, so that, from Theorem 2.4, the function $G$ is strictly concave and $G(0) = 0$. The argument above then gives $G(s) \ge s$ for all $s$, which therefore precludes the existence of any $s > 0$ with $G(s) = s$. This completes the proof. □

REMARK 3.5. From (2.31), we see that

$$(3.14) \quad G'(0) = \gamma \int_0^\infty \mathbb{E}^{(1)} Z^{(0)}_t \, dt.$$

Thus, $G'(0)$ is the average number of successful propagules produced in a patch colonized by a single immigrant, before population extinction in that patch, disregarding other colonizations. This number may be considered a reproduction number for colonizers of an empty habitat, as used in epidemic models [9], thus, $G'(0) > 1$ is the natural condition to ensure (meta)population persistence. Indeed, a similar condition has been presented by Chesson [7] and Casagrandi and Gatto [6]. See also [21], in which an analogous quantity is used as the invasion fitness of a mutant; a discussion along their lines is, however, rather beyond the scope of this paper.

Comparing the process $Z^{(0)}_t$ with a process $\widehat{Z}^{(0)}_t$ which has $\hat{b}_j = b_0$ and $\hat{d}_j = d_0$, one immediately obtains

$$\mathbb{E}^{(1)} Z^{(0)}_t \le e^{-\nu t} e^{(b_0 - d_0 - \gamma)t}.$$

Hence, if

$$b_0 - d_0 - \nu \le 0,$$

one has $G'(0) \le 1$.

**4. Convergence to equilibrium.** In this section we prove the convergence of the solutions of (1.2) to the unique positive equilibrium, when it exists, or otherwise to the extinction equilibrium given by $e^0 = (1, 0, 0, \ldots)^T$. Conditions (H1) and (H2) are assumed to hold throughout the section. We begin with two natural bounds on the mean patch size, the first of which bounds $s(t)$ away from infinity.

LEMMA 4.1. *Let $p^0 \geq 0$ and let*
$$s(t) = \sum_{j=0}^{\infty} j p_j(t).$$
*Then*
$$\limsup_{t \to \infty} s(t) < +\infty.$$

PROOF. Multiplying both sides of (1.2) by $i$ and summing for $i$ from 1 to $\infty$, we obtain

$$(4.1) \qquad s'(t) = \sum_{j=0}^{\infty} j b_j p_j(t) - \sum_{j=0}^{\infty} j d_j p_j(t) - \nu s(t).$$

Note that, as in the previous section, the interchange of derivatives and sums is justified, if $p^0 \in D(A_2)$, by the fact that the solution $p(t) \in C^1([0,T]; m^1)$ and satisfies (3.4). By density, (4.1) then holds for all $p^0 \in D(A)$.

Now, using the concavity of $xb(x)$ and the convexity of $xd(x)$ as in (3.8), we have, from (4.1),

$$(4.2) \qquad s'(t) \leq \{b(s(t)) - d(s(t)) - \nu\} s(t).$$

By standard comparison arguments, we easily obtain
$$\limsup_{t \to \infty} s(t) \leq \tilde{s},$$
where
$$\tilde{s} = \inf\{s > 0 : b(s) < d(s) + \nu\}.$$
The set is not empty because of (H2). □

The next lemma gives the complementary comparison result, bounding $s(t)$ away from 0 when $G'(0) > 1$ and $p^0 \neq e^0$. Its proof is very much more difficult, and is the subject of Section 5.

LEMMA 4.2. *Let $G'(0) > 1$ and $d_\infty < +\infty$. If $p^0 \in C$, $p^0 \neq e_0$, then*
$$\liminf_{t \to \infty} s(t) > 0.$$

Now, if $G'(0) > 1$, let $s^*$ be the unique positive fixed point of $G$, as in Theorem 3.1; if $G'(0) \leq 1$, let $s^* = 0$. In the next two lemmas, we show that $s(t)$ converges to $s^*$.

LEMMA 4.3. *Under the same assumptions as in Lemma 4.1, we have*

$$\limsup_{t \to \infty} s(t) \leq s^*.$$

PROOF. Assume, if possible, that

$$\limsup_{t \to \infty} s(t) = \bar{s} > s^*.$$

From the proof of Theorem 3.1, we then have $G(\bar{s}) < \bar{s}$. Choose $\varepsilon$ such that $G(\bar{s} + \varepsilon) < \bar{s}$, and then choose $t_0$ such that $s(t) \leq \bar{s} + \varepsilon$ for all $t \geq t_0$.

If we take $s(t)$ as a fixed given function, we see that the solution of (1.2) can be interpreted as the distribution of an immigration, birth, death and catastrophe process $Z(t)$ with time varying immigration rate $s(t)$, starting at time $t_0$ with distribution $p(t_0)$. By an easy stochastic comparison (see [3]), that process is dominated in $t \geq t_0$ by a process $Z^{(\bar{s}+\varepsilon)}$ with constant immigration rate $\bar{s} + \varepsilon$ and with the same initial condition $p(t_0)$. In Theorem 2.3, it is shown that $Z^{(\bar{s}+\varepsilon)}$ is positive recurrent and that its equilibrium distribution has finite mean $G(\bar{s} + \varepsilon)$ as in (3.2); furthermore, from Theorem 2.3 and from (2.19) with $\delta = 0$, it follows that

$$\lim_{t \to \infty} \mathbb{E} Z^{(\bar{s}+\varepsilon)}(t) = G(\bar{s} + \varepsilon) \tag{4.3}$$

if $\sum_{j \geq 1} j p_j(t_0) < \infty$, true for all $p^0 \in \mathcal{D}(A)$ because of Theorem A.

Hence, if $p(t)$ is the distribution of $Z(t)$, we have, using also (2.18),

$$\limsup_{t \to \infty} \sum_{j=1}^{\infty} j p_j(t) \leq G(\bar{s} + \varepsilon) < \bar{s}.$$

On the other hand, $p(t)$ is the solution of (1.2) and $s(t)$ was defined as $\sum_{j=0}^{\infty} j p_j(t)$. The previous inequality thus reads

$$\limsup_{t \to \infty} s(t) = \limsup_{t \to \infty} \sum_{j=1}^{\infty} j p_j(t) < \bar{s},$$

contradicting $\bar{s} = \limsup_{t \to \infty} s(t)$. □

The companion result is as follows.

LEMMA 4.4. *Let $G'(0) > 1$ and $d_\infty < +\infty$. If $p^0 \in C$, $p^0 \neq e_0$, then*

$$\liminf_{t \to \infty} s(t) \geq s^*.$$

PROOF. Assume, if possible, that
$$\liminf_{t\to\infty} s(t) = \bar{s} < s^*.$$

Then, from Lemma 4.2, we have $0 < \bar{s} < s^*$. Since $G(\bar{s}) > \bar{s}$, as seen in the proof of Theorem 3.1, we can choose $\varepsilon$ such that $G(\bar{s} - \varepsilon) > \bar{s}$.

As in the proof of Lemma 4.3, choosing $t_0$ such that $s(t) \geq \bar{s} - \varepsilon$ for all $t \geq t_0$, we can compare the process with immigration rate $s(t)$ to the process with immigration rate $\bar{s} - \varepsilon$. In this way, we obtain

$$\liminf_{t\to\infty} s(t) = \liminf_{t\to\infty} \sum_{j=1}^{\infty} j p_j(t) \geq G(\bar{s} - \varepsilon) > \bar{s},$$

reaching a contradiction. □

Combining these lemmas, we can prove the following theorem.

THEOREM 4.5. *Let* (H1) *and* (H2) *hold, and let* $p^0 \in C \setminus \{e^0\}$. *Then the solution of* (1.2) *converges to the unique positive equilibrium, if* $G'(0) > 1$ *and* $d_\infty < +\infty$, *and to* $e^0$ *if* $G'(0) \leq 1$.

PROOF. The previous lemmas together yield
$$\lim_{t\to\infty} s(t) = s^*.$$

Now, the interpretation of $p(t)$ as the distribution at time $t$ of an immigration, birth, death and catastrophe process $Z$ with immigration rate $\gamma s(t)$ shows, as in the proof of Lemma 4.3, that $p(t)$ is asymptotically bounded between the distributions of the processes $Z^{(s^*-\varepsilon)}$ and $Z^{(s^*+\varepsilon)}$ for any $\varepsilon > 0$; that is, for any $l \geq 0$,

$$\sum_{j \geq l} \pi_j^{(s^*-\varepsilon)} \leq \liminf_{t\to\infty} \sum_{j \geq l} p_j(t) \leq \limsup_{t\to\infty} \sum_{j \geq l} p_j(t) \leq \sum_{j \geq l} \pi_j^{(s^*+\varepsilon)}.$$

But Theorem 2.4 implies the continuity in $s$ of $\pi^{(s)}(f)$ with $f = \mathbb{1}_{[l,\infty)}$, proving the theorem. □

REMARK 4.6. The condition $d_\infty < +\infty$ is used in the proof of Lemma 4.2. There is no reason to suppose that it is necessary for Theorem 4.5 to be true, but our proof makes essential use of it.

**5. Repulsion from the extinction equilibrium.** The aim of this section is to prove Lemma 4.2. To do so, we employ a result from the theory of persistence, which we now recall.

Let $X$ be a metric space (with metric $d$) which is the disjoint union of two sets $X_1$ and $X_2$, and suppose that $\Phi$ is a continuous semiflow on $X_1$. Thieme [25] gives the following definitions:

- $X_2$ is a *weak repeller* for $X_1$ if

$$\limsup_{t\to\infty} d(\Phi_t(x_1), X_2) > 0 \qquad \forall\, x_1 \in X_1.$$

- $X_2$ is a *uniform weak repeller* for $X_1$ if there exists $\varepsilon > 0$ such that

$$\limsup_{t\to\infty} d(\Phi_t(x_1), X_2) > \varepsilon \qquad \forall\, x_1 \in X_1.$$

- $X_2$ is a *strong repeller* for $X_1$ if

$$\liminf_{t\to\infty} d(\Phi_t(x_1), X_2) > 0 \qquad \forall\, x_1 \in X_1.$$

- $X_2$ is a *uniform strong repeller* for $X_1$ if there exists $\varepsilon > 0$ such that

$$\liminf_{t\to\infty} d(\Phi_t(x_1), X_2) > \varepsilon \qquad \forall\, x_1 \in X_1.$$

In our application, the space $X$ will be the convex set

$$C = \left\{ p \in m^1 : p \geq 0,\ \sum_{j=0}^{\infty} p_j = 1 \right\}$$

with

$$d(p,q) = |p_0 - q_0| + \sum_{j=1}^{\infty} j |p_j - q_j|,$$

and the continuous semiflow $\Phi_t(p) = \Phi(t,p)$ is given by the solution $p(t)$ of (1.6) with $p(0) = p$. We take $X_2$ to be $\{e^0\}$ and $X_1 := C \setminus \{e^0\}$; with these definitions, the thesis of Lemma 4.2 is that $X_2$ is a strong repeller for $X_1$.

To prove the lemma, we use Theorem 6.2 of [25], which we state in a form *simplified to our present needs.*

THEOREM B ([25]). *Let $X$ be a metric space which is the disjoint union of the two sets $X_1$ (open in $X$) and $X_2$; let $\Phi$ be a continuous semiflow on $X_1$. Assume the following:*

(A) *There exists a subset $Y_1 \subset X_1$ such that, for all $x \in X_1$, there exists $t(x) > 0$ such that $\Phi_t(x) \in Y_1$ for all $t > t(x)$.*
($C_{6,1}$) *For any $y \in Y_1$, the orbit $\Phi([0,\infty) \times \{y\})$ has compact closure.*
($C_{6,2}$) *$\bigcup_{y \in Y_1} \omega(y)$ has compact closure, where, as usual, $\omega(y)$ is the $\omega$-limit set.*
(R) *The set $Y_1 \cap \{x \in X;\ d(x, X_2) = \varepsilon\}$ is bounded.*

*Then $X_2$ is a uniform strong repeller whenever it is a uniform weak repeller.*

We prove that $X_2$ is a uniform weak repeller, and then Theorem B lets us conclude that $X_2$ is a (uniform) strong repeller, which is the thesis of Lemma 4.2. To start with, we show that the assumptions of Theorem B are satisfied.

Lemma 4.1 shows that, if we choose

$$Y_1 = \left\{ y : \sum_i i y_i \leq \tilde{s} \right\},$$

then assumption (A) holds. Indeed, the proof of Lemma 4.1 shows that, if $y \in Y_1$, then $\Phi_t(y) \in Y_1$ for all $t \geq 0$. Assumption (R) is immediate, because $X_2$ is bounded. The following lemma establishes the other two assumptions of Theorem B. For its proof, note that a set $E \subset m^1$ has compact closure if (and only if) $\lim_{N \to \infty} \sum_{n=N}^{\infty} i|x_i| = 0$ uniformly for $x \in E$; that is, if, given any $\varepsilon > 0$, there exists $N = N(\varepsilon) \geq 1$ such that $\sum_{n=N}^{\infty} i|x_i| < \varepsilon$ for all $x \in E$.

LEMMA 5.1. *If the continuous semi-flow $\Phi$ is given by the solutions $p(t)$ of (1.6) and $C$, $X_1$ and $X_2$ are the sets defined above, then assumptions $(C_{6,1})$ and $(C_{6,2})$ hold.*

PROOF. As in Section 4, observe that $p(t) = \Phi_t(y)$ is the distribution of an immigration, birth, death and catastrophe process $Z_t$ with immigration rate $s(t)$ starting at time 0 with distribution $y$. If $y \in Y_1$, this is dominated by an immigration, birth, death and catastrophe process $Z_t^{(\tilde{s})}$ with constant immigration rate $\tilde{s}$ (because of the previous remark), whose transition probabilities we denote by

$$\tilde{p}_{ij}(t) = \mathbb{P}(Z_t^{(\tilde{s})} = j | Z_0^{(\tilde{s})} = i).$$

Stochastic comparison then gives

$$(5.1) \quad \sum_{n=N}^{\infty} n p_n(t) \leq \sum_{n=N}^{\infty} n \sum_{i=0}^{\infty} y_i \tilde{p}_{in}(t) = \sum_{i=0}^{\infty} y_i \sum_{n=N}^{\infty} n \tilde{p}_{in}(t).$$

To estimate the right-hand side, we use (2.15) in Theorem 2.3; choosing $\delta$ such that $c(1+\delta) < \nu$, we obtain

$$\sum_{n=N}^{\infty} n \tilde{p}_{in}(t) \leq \frac{1}{N^\delta} \sum_{n=N}^{\infty} n^{1+\delta} \tilde{p}_{in}(t) \leq \frac{1}{N^\delta} \mathbb{E}^{(i)}(Z_t^{(\tilde{s})})^{1+\delta} \leq \frac{C_i}{N^\delta},$$

uniformly for all $t \geq 0$, where $C_i$ is a constant depending only on $i$. Note also that, for $\delta = 0$, (2.15) implies that

$$(5.2) \quad \sum_{n=0}^{\infty} n \tilde{p}_{jn}(t) = \mathbb{E}^{(j)}(Z_t^{(\tilde{s})}) \leq K_1(j+1).$$

To prove $(C_{6,1})$, take $y \in Y_1$; choose $\varepsilon > 0$. Find $N_1$ such that

$$\sum_{i=N_1}^{\infty} i y_i \leq \frac{\varepsilon}{4K_1},$$

and pick $N_2$ such that $N_2 > (\frac{2C_i}{\varepsilon})^{1/\delta}$ for $i = 0, \ldots, N_1$; then $N_2$ is the required constant. In fact, using (5.2), we obtain

$$\sum_{i=0}^{\infty} y_i \sum_{n=N_2}^{\infty} n\tilde{p}_{in}(t) = \sum_{i=0}^{N_1} y_i \sum_{n=N_2}^{\infty} n\tilde{p}_{in}(t) + \sum_{i=N_1}^{\infty} y_i \sum_{n=N_2}^{\infty} n\tilde{p}_{in}(t)$$

$$\leq \sum_{i=0}^{N_1} y_i \frac{C_i}{N_2^{\delta}} + \sum_{i=N_1}^{\infty} y_i(i+1)K_1 \leq \frac{\varepsilon}{2} \sum_{i=0}^{N_1} y_i + 2K_1 \sum_{i=N_1}^{\infty} iy_i \leq \varepsilon.$$

In order to prove $(C_{6,2})$, we prove that, for any $\varepsilon > 0$, there exists $N = N(\varepsilon) \geq 1$ such that, for all $y \in Y_1$, there exists $t_0 = t_0(y)$ such that

$$\sum_{n=N}^{\infty} n \sum_{i=0}^{\infty} y_i \tilde{p}_{in}(t) < \varepsilon \qquad \text{for all } t \geq t_0.$$

Indeed, assume that this is true, and take $q \in \omega(y)$ for some $y \in Y_1$. Then there exists a sequence $\{t_k\}$ with $t_k \to \infty$ such that

(5.3) $$\sum_{n=0}^{\infty} n|p_n(t_k) - q_n| \to 0 \qquad \text{as } k \to \infty.$$

Take $k$ such that $t_k > t_0(y)$ and that the difference in (5.3) is less than $\varepsilon$. Then

$$\sum_{n=N}^{\infty} nq_n \leq \sum_{n=N}^{\infty} n|p_n(t_k) - q_n| + \sum_{n=N}^{\infty} p_n(t_k)$$

$$\leq \sum_{n=0}^{\infty} n|p_n(t_k) - q_n| + \sum_{n=N}^{\infty} n \sum_{i=0}^{\infty} y_i \tilde{p}_{in}(t) < 2\varepsilon,$$

using also (5.1), so that $(C_{6,2})$ is proved.

Now choose $\delta > 0$ such that $c(1 + \delta) < \nu$, and recall as above that, for each $i$, there exists $C_i < \infty$ such that $\limsup_{t\to\infty} \mathbb{E}^{(i)}(Z_t^{(\tilde{s})})^{1+\delta} \leq C_i$. Hence, for each $i$, there exists $t_0(i)$ such that $\mathbb{E}^{(i)}(Z_t^{(\tilde{s})})^{1+\delta} \leq 2C_i$ for all $t \geq t_0$ and, hence, that

$$\sum_{n=N}^{\infty} n\tilde{p}_{in}(t) \leq 2C_i N^{-\delta} \qquad \text{for all } t \geq t_0(i).$$

Fix $\varepsilon > 0$. Choose $y \in Y_1$ and find $N_1 = N_1(\varepsilon, y)$ such that

$$\sum_{i=N_1}^{\infty} iy_i < \frac{\varepsilon}{4K_1},$$

where $K_1$ is as in (5.2); set $t_0(y) = \max_{i=0,\ldots,N_1} t_0(i)$, and choose

$$N = \left\{ 4\varepsilon^{-1} \max_{1 \leq i \leq N_1} C_i \right\}.$$

Then, for $t \geq t_0(y)$, we have

$$\sum_{i=0}^{\infty} y_i \sum_{n=N}^{\infty} n\tilde{p}_{in}(t) = \sum_{i=0}^{N_1} y_i \sum_{n=N}^{\infty} n\tilde{p}_{in}(t) + \sum_{i=N_1+1}^{\infty} y_i \sum_{n=N}^{\infty} n\tilde{p}_{in}(t)$$

$$\leq \sum_{i=0}^{N_1} y_i \frac{\varepsilon}{2} + \sum_{i=N_1+1}^{\infty} y_i(i+1)K_1 \leq \varepsilon,$$

proving (5.3). □

Now we prove that $\{e^0\}$ is a weak repeller through linearization. Since we restrict our considerations to vectors $p(t)$ in the convex set $C$, we have $p_0(t) = 1 - \sum_{j=1}^{\infty} p_j(t)$. Hence, we need only examine the vector $(p_1, p_2, \ldots)^T$. With a slight abuse of notation, we now set

$$X = \left\{ x = (x_1, x_2, \ldots)^T \in \ell^1 : \sum_{j=1}^{\infty} j|x_j| < +\infty \right\}$$

with norm $\|x\| = \sum_{j=1}^{\infty} j|x_j|$, noting that $e^0$ now translates into the point 0 of $X$, and we denote here by $A$ and $F$ the operators defined in (1.4)–(1.5) but restricted to $X$, and using $p_0 = 1 - \sum_{j=1}^{\infty} p_j$ in the definition of $F$. We then define $X_+$ to be the nonnegative cone in $X$; note that $X_+$ is the counterpart of the convex set $C$ defined above.

Equation

(5.4) $$p' = Ap + F(p)$$

corresponding to (1.6) now has 0 as the equilibrium, corresponding to the extinction equilibrium $e^0$ of (1.2). We again use $\Phi_t(u^0)$ to denote the solution of (5.4) satisfying $u(0) = u^0$. This corresponds to the semi-flow of Lemma 5.1, except that we now neglect the 0th component. Note that the metric in the convex set $C$ is equivalent to the norm in $X$, since

$$d(u, v) = |u_0 - v_0| + \sum_{i=1}^{\infty} i|u_i - v_i|$$

$$= \left| \left(1 - \sum_{i=1}^{\infty} u_i\right) - \left(1 - \sum_{i=1}^{\infty} v_i\right) \right| + \sum_{i=1}^{\infty} i|u_i - v_i|$$

$$\leq 2 \sum_{i=1}^{\infty} i|u_i - v_i| = 2\|u - v\|_X,$$

while obviously $\|u - v\|_X \leq d(u, v)$.

Note also that $A$ is the generator of a defective Markov process, the process $Z_t^{(0)}$ of Section 2 restricted to the state space $\mathbb{N} \setminus \{0\}$. In the rest of this section, we only consider processes with zero immigration rate; thus, when there is no ambiguity, we drop the superscript $^{(0)}$ and denote by $Z_t$ the process with zero immigration rate.

From the results of Section 2, one immediately sees that $Z_t$ is exponentially absorbed at 0; more precisely, (2.13) with $f = e$ and $l = s = 0$ implies, in the present notation, that

$$\|e^{At}\| \leq Ce^{-\alpha t} \tag{5.5}$$

for some positive constants $C$ and $\alpha$. This implies that

$$\{\mathbb{R}\lambda > -\alpha\} \subset \rho(A);$$

moreover, we have the representation

$$((\lambda - A)^{-1}v)_i = \sum_j v_j \hat{P}_{ji}(\lambda), \tag{5.6}$$

where "ˆ" denotes the Laplace transform and $P_{ji}(t)$ is $\mathbb{P}(Z_t = i | Z_0 = j)$.

We now discuss the stability of the 0 equilibrium of (5.4) using the linearization principle. We first note that

$$F'(0)u = \varphi(u)e^1 \qquad \forall u \in X, \tag{5.7}$$

where

$$\varphi(u) = \gamma \sum_j j u_j \tag{5.8}$$

and $e^1 = (1, 0, 0, \dots)^T$. Since $F'(0)$ is one-dimensional, hence, compact, the essential spectrum [27] of $A + F'(0)$ coincides with that of $A$, which, from (5.5), is less or equal than $-\alpha$. The type of the semigroup $e^{(A+F'(0))t}$ can then be understood from the spectrum of $A + F'(0)$.

Using (5.7), we can establish, through direct computation, the following result.

LEMMA 5.2. *If $\lambda$ is in $\rho(A)$, then $\lambda$ belongs to $\rho(A + F'(0))$ if and only if $\varphi((\lambda - A)^{-1}e^1) \neq 1$. In that case,*

$$(\lambda - A - F'(0))^{-1}v = (\lambda - A)^{-1}v + \frac{\varphi((\lambda - A)^{-1}v)}{1 - \varphi((\lambda - A)^{-1}e^1)}(\lambda - A)^{-1}e^1. \tag{5.9}$$

*On the other hand, if*

$$\varphi((\lambda - A)^{-1}e^1) = 1, \tag{5.10}$$

*then $\lambda$ is an eigenvalue with corresponding eigenvector $v = (\lambda - A)^{-1}e^1$.*

From this lemma, we see that an important role is played by the roots of (5.10) in the half-plane $\{\mathbb{R}\lambda > -\alpha\}$. Using the representation (5.6) and standard results on the Laplace transform, as used, for instance, in the theory of age-dependent populations [16], we have the following:

LEMMA 5.3. *There exists, at most, one real root $\lambda_0 > -\alpha$ of (5.10). If $\lambda_0$ exists, all the other roots $\lambda$ satisfy $\mathbb{R}\lambda < \lambda_0$; if there is no real root, there are no complex roots in $\{\mathbb{R}\lambda > -\alpha\}$. In any strip $\{a \leq \mathbb{R}\lambda \leq b\}$, there are, at most, finitely many roots.*

*Finally, if*

$$R_0 := \sum_i i\hat{P}_{1i}(0) = \sum_i i \int_0^\infty P_{1i}(t)\,dt > [=] 1,$$

*then $\lambda_0 > [=] 0$; on the other hand, if $R_0 < 1$, if there is a real root $\lambda_0$, it satisfies $\lambda_0 < 0$.*

REMARK 5.4. Note that

$$R_0 = \sum_i i \int_0^\infty P_{1i}(t)\,dt = G'(0),$$

with $G$ as given in (3.2).

From here on we assume that $R_0 > 1$. Hence, the real eigenvalue $\lambda_0$ is positive. We denote by $\lambda_1, \ldots, \lambda_k$ (with $k \geq 0$) the other roots of (5.10) such that $\mathbb{R}\lambda_j > 0$, and by $\lambda_{k+1}, \ldots, \lambda_n$ (with $n \geq k$) the roots such that $\mathbb{R}\lambda_j = 0$. Since the continuous spectrum (if it exists) of $A + F'(0)$ is contained in $\{\mathbb{R}\lambda \leq -\alpha\}$, we can split the spectrum of $A + F'(0)$ in three spectral sets $\sigma^u = \{\lambda_0, \lambda_1, \ldots, \lambda_k\}$, $\sigma^c = \{\lambda_{k+1}, \ldots, \lambda_n\}$ and $\sigma^s = \{\lambda \in \sigma(A + F'(0)) : \mathbb{R}\lambda < 0\}$.

By standard results (see Theorem III.6.17 in [17]), $X$ can be split into the direct sum of three subspaces $X^u$, $X^c$ and $X^s$, all invariant under $A + F'(0)$. Moreover, $X^u$ and $X^c$ are finite-dimensional ($X^u$ includes at least $v_0$, the eigenvector corresponding to $\lambda_0$, while $X^c$ may well consist only of 0). This would be enough to establish the instability of the 0 equilibrium. However, we wish to prove that all initial data $u \geq 0$, $u \neq 0$, are repelled away from 0, and this requires further work.

The following lemma uses the results of Bates and Jones [4] to establish the existence of unstable and centre stable manifolds $W^u$ and $W^{cs}$ for equation (5.4) at 0. The conditions of their Theorem 1.2 are satisfied in view of Arrigoni's results, as summarized at the end of Section 1, together with (5.5) and the properties of the eigenspaces discussed following Lemma 5.3.

Defining $X^{cs} = X^c \oplus X^s$, and letting $P^u$ and $P^{cs}$ denote the corresponding projections, [4], Theorem 1.2 and its consequence (P3) yield the following result.

LEMMA 5.5. *There exist a neighborhood $U \ni 0$ and Lipschitz functions $h^u \colon P^u(U) \to X^{cs}$ and $h^{cs} \colon P^{cs}(U) \to X^u$ with $h^u(0) = (h^u)'(0) = h^{cs}(0) = (h^{cs})'(0) = 0$ such that*

$$W^u = \{u^u + h^u(u^u) \colon u^u \in P^u(U)\}$$

*is the unstable manifold (in $U$) of 0, and*

$$W^{cs} = \{u^{cs} + h^{cs}(u^{cs}) \colon u^{cs} \in P^{cs}(U)\}$$

*is a centre-stable manifold.*

*Furthermore, there exists a neighborhood $V \subset U$ of 0 such that, if $u^0 \in V \setminus W^{cs}$, then there exists $\tau > 0$ such that $\Phi_\tau(u^0) \notin V$.*

The final statement of the lemma shows that, if a solution comes close enough to 0 to be in the neighborhood $V$, and if it is then at a point not in $W^{cs}$, then it has to leave $V$ at some later time. Hence, the limes superior of any solution curve is necessarily positive, if it can be established that, for some $\varepsilon > 0$, *no* points of $X_+ \cap B_\varepsilon$ except for 0 are in $W^{cs}$, where $B_\varepsilon$ denotes the ball of radius $\varepsilon$ centred at 0. If this is the case, then $\{0\}$ is a uniform weak repeller for $X_+ \setminus \{0\}$ in the system (5.4), which is equivalent to $\{e^0\}$ being a weak repeller for $C \setminus \{e^0\}$ in (1.6). Applying Theorem B, Lemma 4.2 would then follow.

To show that indeed $W^{cs} \cap X_+ \cap B_\varepsilon = \{0\}$ for some $\varepsilon > 0$, we begin by writing the eigenprojections explicitly.

LEMMA 5.6. *The projection $P_0$ on the eigenspace corresponding to $\lambda_0$ is*

$$P_0 v = -\frac{\varphi((\lambda_0 - A)^{-1} v)}{\varphi'((\lambda_0 - A)^{-1} e^1)} (\lambda_0 - A)^{-1} e^1.$$

*The projection $P^u$ on $X^u$ is given by*

$$P^u v = P_0 v + \sum_{j=1}^{k} \frac{1}{2\pi} \oint_{\Gamma_j} \frac{\varphi((\lambda - A)^{-1} v)}{1 - \varphi((\lambda - A)^{-1} e^1)} (\lambda - A)^{-1} e^1 \, d\lambda,$$

*where $\Gamma_j$ is a circle around $\lambda_j$ that does not include other elements of the spectrum.*

PROOF. It follows from the construction of the projection operators as in formula (III.6.19) of [17] and from (5.9). □

Note that

$$\varphi'((\lambda_0 - A)^{-1} e^1) = -\int_0^\infty t e^{-\lambda_0 t} \sum_{i=1}^\infty i P_{1i}(t) \, dt < 0.$$

On the other hand, it may well be $\varphi'((\lambda_j - A)^{-1}e^1) = 0$ when $1 \leq j \leq n$, so that the other projections may have a more complex form.

As a consequence, we immediately have the following result.

LEMMA 5.7. *If $v \in X^{cs}$, then $\varphi((\lambda_0 - A)^{-1}v) = 0$.*

PROOF. The explicit representation of $P_0$ shows that if $\varphi((\lambda_0 - A)^{-1}v) \neq 0$, then $P_0 v \neq 0$. However, $v \in X^{cs}$ implies that $P_0 v = 0$. □

This lemma implies that $X^{cs} \cap X_+ = \{0\}$, because if $v \geq 0$ and $v \neq 0$, then $(\lambda_0 - A)^{-1}v > 0$ by (5.6); hence, it follows from (5.8) that $\varphi((\lambda_0 - A)^{-1}v) > 0$ also. This is almost what we need, since $W^{cs}$ and $X^{cs}$ are close to one another near 0, and we are thus close to showing that $W^{cs} \cap X_+ \cap B_\varepsilon = \{0\}$ for some $\varepsilon > 0$. To make the transition from $X^{cs}$ to $W^{cs}$, we first show that, for $v \geq 0$, $\varphi((\lambda_0 - A)^{-1}v)$ is large enough.

LEMMA 5.8. *Assume that $d_\infty < \infty$, and take $v \geq 0$. Then*

$$(5.11) \qquad \varphi((\lambda_0 - A)^{-1}v) \geq \frac{\|v\|}{d_\infty + \gamma + \nu + \lambda_0}.$$

PROOF. We start from the identity

$$\sum_i i \hat{P}_{ji}(\lambda_0) = \int_0^\infty e^{-\lambda_0 t} \mathbb{E}^{(j)}(Z_t)\,dt.$$

An easy coupling argument shows that $Z_t$ is stochastically larger than a death-and-catastrophe process with death rate $d_\infty + \gamma$. Hence,

$$\int_0^\infty e^{-\lambda_0 t} \mathbb{E}_j(Z_t)\,dt \geq j \int_0^\infty e^{-(\lambda_0 + \gamma + d_\infty + \nu)t}\,dt = \frac{j}{\lambda_0 + \gamma + d_\infty + \nu}.$$

Now, if $v \geq 0$, we have

$$\varphi((\lambda_0 - A)^{-1}v) = \sum_{i,j} i v_j \hat{P}_{ji}(\lambda_0)$$

$$\geq \sum_j v_j \frac{j}{\lambda_0 + \gamma + d_\infty + \nu} = \frac{\|v\|}{d_\infty + \gamma + \nu + \lambda_0}. \qquad \square$$

Using the above lemma together with Lemma 5.7, we can now show that the norm of $v^-$ is quite large, whenever $v \in X^{cs}$. Here, $v^-$ denotes the negative part of $v$: $v = v^+ - v^-$, with $(v^+)_i = \max\{0, v_i\}$ and $(v^-)_i = \max\{0, -v_i\}$.

LEMMA 5.9. *Assume that $d_\infty < \infty$. If $v \in X^{cs}$, there exists $\eta > 0$ such that $\|v^-\| \geq \eta \|v\|$.*

PROOF. From Lemma 5.7, we have

$$0 = \varphi((\lambda_0 - A)^{-1}v) = \varphi((\lambda_0 - A)^{-1}v^+) - \varphi((\lambda_0 - A)^{-1}v^-)$$

$$\geq \frac{\|v^+\|}{d_\infty + \gamma + \nu + \lambda_0} - \|(\lambda_0 - A)^{-1}\| \, \|v^-\|,$$

using Lemma 5.8 and the obvious identity $\|\varphi\| = 1$. Hence, using $\|v^+\| = \|v\| - \|v^-\|$,

$$\left(\|(\lambda_0 - A)^{-1}\| + \frac{1}{d_\infty + \gamma + \nu + \lambda_0}\right)\|v^-\| \geq \frac{1}{d_\infty + \gamma + \nu + \lambda_0}\|v\|,$$

which yields the thesis. $\square$

We now use this result, together with the closeness of $X^{cs}$ and $W^{cs}$, to conclude that $W^{cs} \cap X_+ \cap B_\varepsilon = \{0\}$ for some $\varepsilon > 0$.

LEMMA 5.10. *Assume that $d_\infty < \infty$. Then there exists $\varepsilon > 0$ such that $v \geq 0$, $v \in B_\varepsilon \cap W^{cs}$ implies $v = 0$.*

PROOF. First take $\delta$ such that $\|v^{cs}\| \leq \delta$ implies $\|h^{cs}(v^{cs})\| \leq \frac{\eta}{2}\|v^{cs}\|$. Then take $\varepsilon = \delta/\|P^{cs}\|$. Assume that $v = v^{cs} + h^{cs}(v^{cs}) \geq 0$ with $\|v\| \leq \varepsilon$. Then it follows that $\|v^{cs}\| = \|P^{cs}(v)\| \leq \delta$.

Split $v^{cs} = (v^{cs})^+ - (v^{cs})^-$. Then we have

$$\sum_{i:\, v_i^{cs} < 0} iv_i = \sum_{i:\, v_i^{cs} < 0} i[v_i^{cs} + (h^{cs}(v^{cs}))_i] \leq -\sum_{i=1}^\infty i(v^{cs})_i^- + \sum_{i=1}^\infty i|(h^{cs}(v^{cs}))_i|$$

$$= -\|(v^{cs})^-\| + \|h^{cs}(v^{cs})\| \leq -\eta\|v^{cs}\| + \frac{\eta}{2}\|v^{cs}\|,$$

using Lemma 5.8 and $\|v^{cs}\| \leq \delta$. This contradicts with $v \geq 0$ unless $v^{cs} = v = 0$.

$\square$

We have now proved what we need to show that $\{0\}$ is a uniform weak repeller for $X_+ \setminus \{0\}$. The details are as follows. We recall that we have $R_0 = G'(0) > 1$.

LEMMA 5.11. *Assume that $d_\infty < \infty$. Then there exists $\varepsilon_0$ such that for all $u^0 \geq 0$, $u^0 \neq 0$,*

$$\limsup_{t \to \infty} \|\Phi_t(u^0)\| \geq \varepsilon_0.$$

PROOF. Take

$$\varepsilon_0 = \tfrac{1}{2}\min\biggl\{\varepsilon, \inf_{v\in X\setminus V}\|v\|\biggr\},$$

where $\varepsilon$ is that of Lemma 5.10, while $V$ is that of Lemma 5.5.

Assume that $\|\Phi_t(u^0)\| < 2\varepsilon_0$ for all $t \geq t_0$. Since $u^0 \geq 0$, the invariance of the positive cone under (5.4) gives $\Phi_t(u^0) \geq 0$; moreover, $\Phi_t(u^0) \neq 0$. Hence, Lemma 5.10 implies that $\Phi_{t_0}(u^0) \notin W^{cs}$. From Lemma 5.5, it then follows that $\Phi_\tau(u^0) \notin V$ for some $\tau > t_0$, contradicting $\|\Phi_t(u^0)\| < 2\varepsilon_0$ for all $t \geq t_0$. □

PROOF OF LEMMA 4.2. Going back to the semi-flow $\Phi_t$ on $C$, note that

$$d(\Phi_t(u_0), e^0) = |1 - p_0(t)| + \sum_{j=1}^\infty i|p_i(t)| = \sum_{j=1}^\infty p_i(t) + \sum_{j=1}^\infty ip_i(t)$$

$$\leq 2\sum_{j=1}^\infty ip_i(t) = \|\Phi_t(u_0)\|,$$

while obviously

$$\|\Phi_t(u_0)\| \leq d(\Phi_t(u_0), e^0).$$

Hence, Lemma 5.11 states that $\{e^0\}$ is a uniform weak repeller for $C \setminus \{e^0\}$. But now Theorem B, together with Lemma 5.1, yields Lemma 4.2. □

**Acknowledgment.** We are grateful to the referees for suggestions which have led to a much improved presentation.

ABTEILUNG ANGEWANDTE MATHEMATIK
UNIVERSITÄT ZÜRICH
WINTERTHURERSTRASSE 190
CH-8057 ZÜRICH
SWITZERLAND
E-MAIL: a.d.barbour@math.unizh.ch

DIPARTIMENTO DI MATEMATICA
UNIVERSITÀ DI TRENTO
VIA SOMMARIVE 14
38050 POVO (TN)
ITALY
E-MAIL: pugliese@science.unitn.it